\documentclass[10pt,reqno,a4paper]{article}

\usepackage[blocks]{authblk}

\usepackage{amsmath,amsthm,amstext,amscd,amssymb,euscript,url,scalerel,relsize,dsfont,mathrsfs}
\usepackage{mathtools}
\usepackage{graphicx}
\usepackage{subfig}
\usepackage{color}

\newcommand{\R}{\mathds R}
\newcommand{\N}{\mathds N}

\newcommand{\E}{\mathds E}
\newcommand{\Zd}{\mathds Z^d}

\newcommand{\loc}{\scriptscriptstyle{\mathrm{loc}}}
\newcommand{\qloc}{\scriptscriptstyle{\mathrm{qloc}}}
\newcommand{\Loc}{\mathcal{L}(\Omega)}
\newcommand{\Qloc}{\mathcal{QL}(\Omega)}

\renewcommand{\phi}{\varphi}
\newcommand{\epsi}{\ensuremath{\varepsilon}}
\newcommand{\la}{\ensuremath{\Lambda}}
\newcommand{\si}{\ensuremath{\sigma}}
\newcommand{\pee}{\ensuremath{\mathbb{P}}}

\newcommand{\db}{\bar{d}}

\newcommand{\re}{\scaleto{s}{5pt}}

\def\1{\mathds{1}}

\makeatletter
\newcommand{\triplenorm}{\@ifstar\@opnorms\@opnorm}
\newcommand{\@opnorms}[1]{%
  \left|\mkern-1.5mu\left|\mkern-1.5mu\left|
   #1
  \right|\mkern-1.5mu\right|\mkern-1.5mu\right|
}
\newcommand{\@opnorm}[2][]{%
  \mathopen{#1|\mkern-1.5mu#1|\mkern-1.5mu#1|}
  #2
  \mathclose{#1|\mkern-1.5mu#1|\mkern-1.5mu#1|}
}
\makeatother

\usepackage{twoopt}
\newcommandtwoopt{\myfrac}[4][0pt][0pt]{\genfrac{}{}{}{}{\raisebox{#1}{$#3$}}{\raisebox{-#2}{$#4$}}}

\newtheorem{theorem}{{\small T}{\scriptsize HEOREM}}[section]
\newtheorem{corollary}{{\bf{\small C}{\scriptsize OROLLARY}}}[section]
\newtheorem{proposition}{{\bf{\small P}{\scriptsize ROPOSITION}}}[section]
\newtheorem{lemma}{{\bf{\small L}{\scriptsize EMMA}}}[section]
\newtheorem{remark}{{\bf{\small R}{\scriptsize EMARK}}}[section]
\newtheorem{definition}{{\bf{\small D}{\scriptsize EFINITION}}}[section]
\newtheorem{induction}{{\bf{\small I}{\scriptsize NDUCTIVE HYPOTHESIS}}}[section]

\renewenvironment{proof}[1]
{\noindent{{\bf{\small{P}{\scriptsize ROOF}}}.}\hspace{0.1cm} #1} {$\;\qed$\newline}

\newenvironment{sproof}[2]
{\noindent{{\bf{\small{P}{\scriptsize ROOF}}} #1.}\hspace{0.1cm} #2} {$\;\qed$\newline}

\newcommand{\beq}{\begin{eqnarray}}
\newcommand{\eeq}{\end{eqnarray}}

\newcommand{\ba}{\begin{align*}}
\newcommand{\ea}{\end{align*}}

\newcommand{\be}{\begin{equation}}
\newcommand{\ee}{\end{equation}}

\newcommand{\bl}{\begin{lemma}}
\newcommand{\el}{\end{lemma}}

\newcommand{\br}{\begin{remark}}
\newcommand{\er}{\end{remark}}

\newcommand{\bt}{\begin{theorem}}
\newcommand{\et}{\end{theorem}}

\newcommand{\bd}{\begin{definition}}
\newcommand{\ed}{\end{definition}}

\newcommand{\bind}{\begin{induction}}
\newcommand{\eind}{\end{induction}}

\newcommand{\bp}{\begin{proposition}}
\newcommand{\ep}{\end{proposition}}

\newcommand{\bc}{\begin{corollary}}
\newcommand{\ec}{\end{corollary}}

\newcommand{\bi}{\begin{itemize}}
\newcommand{\ei}{\end{itemize}}

\newcommand{\ben}{\begin{enumerate}}
\newcommand{\een}{\end{enumerate}}

\newcommand{\bpr}{\begin{proof}}
\newcommand{\epr}{\end{proof}}

\newcommand{\boS}{\EuScript{S}}

\newcommand{\boB}{\mathcal{B}}

\newcommand{\caC}{{\mathscr C}}

\newcommand{\caG}{{\mathcal G}}

\newcommand{\caK}{{\mathscr K}}
\newcommand{\caL}{{\mathcal L}}

\newcommand{\caP}{{\mathcal P}}

\newcommand{\caT}{{\mathcal T}}
\newcommand{\caU}{{\mathcal U}}

\newcommand{\caX}{{\mathcal X}}

\newcommand{\boM}{{\EuScript M}}

\newcommand{\frakF}{\mathfrak{B}}
\newcommand\Cyl{\mathfrak{C}}

\newcommand{\gcb}[1]{\mathrm{GCB}\!\left(#1\right)}

\newcommand{\dd}{\mathop{}\!\mathrm{d}}
\DeclareMathOperator{\e}{\mathrm{e}}

\usepackage{hyperref}

\definecolor{unbleu}{rgb}{0.03, 0.15, 0.4}

\begin{document}

\title{{\bf Relative entropy, Gaussian concentration and uniqueness of equilibrium states}}

\author[1]{J.-R. Chazottes
\thanks{Email: \texttt{jeanrene@cpht.polytechnique.fr}}}

\author[2]{F. Redig
\thanks{Email: \texttt{F.H.J.Redig@tudelft.nl}}
}

\affil[1]{{\small Centre de Physique Th\'eorique, CNRS, Institut Polytechnique de Paris, France}}
\affil[2]{{\small Institute of Applied Mathematics, Delft University of Technology, The Netherlands}}

\date{Dated: \today}

\maketitle

\begin{abstract}
For a general class of lattice-spin systems, we prove that an abstract Gaussian concentration bound implies positivity of lower entropy density.
As a consequence we obtain uniqueness of translation-invariant Gibbs measures from the Gaussian concentration bound in this general setting. This
extends earlier results from \cite{moles} with a different and very short proof.
\end{abstract}



\section{Introduction}

In \cite{moles} we studied the relation between the Gaussian concentration bound and the uniqueness of equilibrium
states (translation-invariant Gibbs measures) in the context of spin systems on the lattice $\Zd$, where the spin at each lattice site takes a finite
number of values. In this note we prove in a much more general setting that the Gaussian concentration bound implies the strict positivity of the lower relative entropy density. More 
precisely, we obtain a lower bound for lower relative entropy density in terms of a natural distance between translation-invariant probability measures, reminiscent and in the spirit of the 
result of Bobkov and G\"otze \cite{bg} who proved (in a different setting) a lower bound for the relative entropy in terms
of the square of the Wasserstein distance.
Because we work in the thermodynamic limit and are interested in translation-invariant probability measures, there is no translation-invariant distance on the configuration space for which we can apply Bobkov-G\"otze
theorem. We can avoid this problem by introducing a suitable distance on translation-invariant probability measures (rather than on configurations).

We start by proving the lower bound on the lower relative entropy density in the context of general lattice spin systems with state
space $\Omega=S^{\Zd}$ where
the single spins take values in a metric space $S$ of bounded diameter.
The bounded diameter property allows to associate to a quasi-local function $f$ a natural sequence of oscillations $\delta_i f, i\in\Zd$, where $\delta_i f$ represents
the maximal influence on the function $f$ of a change of the spin at site $i$.
In the final section of this paper, we provide a generalization of this by allowing more abstract single spin spaces, and more general
associated sequences of oscillations $\delta_i f, i\in\Zd$.

The rest of our paper is organized as follows.
In Section 2 we introduce the basic setting of lattice spin systems and important function spaces. In Section 3 we introduce the Gaussian concentration bound, relative entropy (density) and formulate and prove our main result in the context of a single spin space with finite diameter. We discuss applications of our result to zero information distance sets, including e.g. the set of equilibrium states w.r.t.
absolutely summable translation-invariant potential.
In Section 4 we consider a generalization introducing an abstract sequence of oscillations, associated Gaussian concentration bound and state and proof the analogue of our main result in this generalized context.
In Section 5 we give various examples.

\section{Setting}


\subsection{Configuration space and the translation operator}\label{buirbu}
We start from a standard Borel space $(S,\mathfrak{b})$ with metric $d_{\scaleto{S}{5pt}}$,
and we let $\mathrm{diam}(S)=\sup_{s,s'\in S} d_S(s,s')$.\footnote{A measurable space $(S,\mathfrak{b})$ is said to be standard Borel if there exists a metric on $S$ which
makes it a complete separable metric space, and $\mathfrak{b}$ then denotes the associated Borel $\sigma$-algebra.}
In the sequel, for notational convenience we omit the symbol $\mathfrak{b}$ and call $S$ a standard Borel space,
where we always assume that the associated $\si$-algebra is the Borel $\si$-algebra $\mathfrak{b}$.

We assume that $\mathrm{diam}(S)<\infty$. Later on, in Section \ref{gensection} we will show how to weaken this assumption.

This space $S$ represents the ``single spin space'', i.e., we will consider lattice-spin configurations where individual spins take values
in $S$.
We denote by $(\Omega_\la,\mathfrak{b}_\la)$ the product space $(S^\la,\prod_{i\in\la} \mathfrak{b}_i\big)$, and $(\Omega,\frakF)$ stands for the lattice-spin configuration space $(\Omega_{\Zd},\mathfrak{b}_{\Zd})$.
Elements of $\Omega$ are called configurations. For $\eta\in\Omega$, we denote by $\eta_i\in S$ its evaluation at site $i\in\Zd$.
By $\sigma_\la$ we mean an element of $\Omega_\la$, and by $\eta_\la\xi_{\la^c}$ a configuration coinciding with $\eta$ on $\la$ and
with $\xi$ on ${\la^c}$. We denote by $\boS$ the set of finite subsets of $\Zd$.

We denote by $\tau_i:\Zd\to\Zd$, $i\in\Zd$, the map which shifts, or translates, by $i$, that is, $\tau_i(j)=j+i$, $j\in\Zd$, and write
$\tau_i(\la)=\la+i=\{j\in\Zd: j=k+i, k\in\Zd\}$, $\la\subset\Zd$.
We define the translation operator  acting on configurations as follows (and use the same symbol). For each $i\in\Zd$, $(\tau_i\si)_j=\si_{j-i}$ for all $j\in\Zd$.
This corresponds to translating $\si$ forward by $i$. We denote by the same symbol the translation operator acting on a function $f:\Omega\to\R$. For each $i\in\Zd$,
$\tau_i f$ is the function defined as $\tau_i f(\si)=f(\tau_i \si)$. A (Borel) probability measure on $\Omega$ is translation invariant if for all $B\in\frakF$ and for all $i\in\Zd$ we have $\mu(\tau_i B)=\mu(B)$.

We denote by $\caP_{\tau}(\Omega)$ the set of translation-invariant probability measures on $\Omega$.

We denote by $\caC(\Omega), \caC_b(\Omega)$ the space of continuous, respectively bounded continuous, real-valued functions.

\subsection{Local oscillations and function spaces}

To a continuous function $f:\Omega\to\R$ we associate a ``sequence'' of ``local oscillations'', $\delta f:=(\delta_i f)_{i\in\Zd}$, defined via
\be\label{difi}
\delta_i f=\sup_{\substack{\si, \eta\,\in\,\Omega: \\ \si_j=\eta_j, \forall j\not= i}} \big(f(\si)-f(\eta)\big).
\ee

Later on, in  Section \ref{gensection}, where we consider the case where $S$ is allowed to have infinite diameter, we will consider a more abstract definition
of $\delta f$. In the case where  $S$ has finite diameter, \eqref{difi} is the most natural choice.

For an integer $p\geq 1$, we define the usual $\ell^p$-norm of $\delta f$,
$
\|\delta f\|^p_p= \sum_{i\in\Zd} (\delta_i f)^p,
$
and $\|\delta f\|_\infty= \sup_i \delta_i f$.

We call a continuous function local if $\delta_i f\not= 0$ for finitely many $i\in D_f\subset\Zd$. $D_f$  is then called the dependence set of $f$.
We denote by $\Loc$ the set of local continuous functions.

We call a continuous function quasi-local if it is the uniform limit of
a sequence of local functions.
If $S$ is compact, then, by the Stone-Weierstrass theorem,  local functions are uniformly dense in $\caC(\Omega)$ and as a consequence quasi-locality and continuity are equivalent.

We denote by $\Qloc$ the space of all continuous quasi-local functions.

For $1\leq p\leq \infty$ the introduce the spaces
\[
\Delta^p(\Omega)= \{ f\in \caC(\Omega): \|\delta f\|_p<\infty\}.
\]

\bl\label{boundlem}
If $f\in \Qloc\cap \Delta^1(\Omega)$ then $f$ is bounded.
If $f\in \Loc\cap\Delta^p(\Omega)$ then $f$ is bounded.
\el
\bpr
Choose $\eta, \si\in \Omega$. We have for every $\la\in\boS$
\be\label{fifof}
|f(\eta)-f(\si)|\leq  |f(\eta)-f(\eta_\la\si_{\la^c})| + \sum_{i\in \la} \delta_i f
\ee
which, upon taking the limit $\la\uparrow\Zd$, using the assumed quasilocality of $f$ gives
\be\label{kokip}
|f(\eta)-f(\si)|\leq \sum_{i\in \Zd} \delta_i f = \|\delta f\|_1<\infty.
\ee
If $f$ is local then we still have the inequality \eqref{fifof} for $\la$ containing the dependence set of $f$.
Because by assumption $\sum_{i} (\delta_i f)^p$ is finite, it follows that $\delta_i f<\infty$ for all $i\in \la$ and therefore $\|\delta f\|_1<\infty$.
Then we obtain \eqref{kokip}, which implies that $f$ is bounded.
\epr

In what follows we will assume that the set of bounded quasi-local functions is measure separating, i.e.,
for two probability measures $\mu\not=\mu'$ we have
that there exists a bounded quasi-local $f$ such that
\[
\int f \dd \mu\neq\int f \dd\mu'.
\]
Because by definition bounded quasi-local functions can be uniformly approximated by bounded local functions, if the set of bounded quasi-local functions is measure separating
then also the set of bounded local functions is measure separating, i.e., for two probability measures $\mu\not=\mu'$ we have
that there exists a bounded local $f$ (which is not constant) such that
\[
\int f \dd \mu\not=\int f \dd\mu'.
\]
We say that $\mu_n\to\mu$ if for all bounded local functions we have $\int f \dd\mu_n\to\int f \dd\mu$ (then by definition of quasilocality, the same holds for bounded quasi-local functions). This induces the so-called 
weak quasi-local topology on probability measures. Notice that in our setting where by assumption the single spin space $S$ is complete metric and separable, this topology coincides with the ordinary weak 
topology, see \cite[p. 898]{enter}.


\section{Gaussian concentration bound and relative entropy}

\subsection{Abstract Gaussian concentration bound}

We can now give the definition of the Gaussian concentration bound in our setting.

\bd
Let $\Omega=S^{\Zd}$ where $S$ is a standard Borel space with finite diameter.
Let $\mu$ be a probability measure on $\Omega$.
We say that $\mu$ satisfies the Gaussian concentration bound with constant $C>0$, abbreviated $\gcb{C}$, if for all bounded local functions $f$
we have 
\be\label{GCB}
\int \e^{f-\int f \dd\mu}\dd\mu\leq \e^{\tfrac{C}{2} \|\delta f\|^2_2}.
\ee
\ed
\br
\leavevmode\\
\textup{(a)}
Observe that the bound \eqref{GCB} does not change if $f$ is replaced by $f+c$, where $c\in\R$ is arbitrary, since $\delta_i(f+c)=\delta_i(f)$ for any $i\in\Zd$. This ``insensitivity'' to constant offsets
in the left-hand side is ensured by the fact that we center $f$ around its expected value. Also observe that \eqref{GCB} is trivially true for functions which are constant.
\newline
\textup{(b)}
We have $\delta_i (\beta f)=|\beta|\, \delta_i f$, for all $i\in\Zd$ and $\beta\in\R$, we thus have
\[
\log\int \e^{\beta(f-\int f \dd\mu)}\dd\mu\leq {\tfrac{C\beta^2}{2} \|\delta f\|^2_2},\;\forall \beta\in\R.
\]
This quadratic upper bound will be crucial in the sequel.

\er
The following proposition shows that \eqref{GCB} automatically extends to a wider class of functions. 
\bp[Self-enhancement of GCB]
Suppose that \eqref{GCB} holds for all bounded local $f$. Then it holds for all $f\in \Qloc\cap \Delta^2(\Omega)$.
\ep
\begin{proof}
By assumption, for $\xi\in\Omega$ fixed,  $f\in \Qloc\cap \Delta^2(\Omega)$ can be uniformly approximated by the local functions
\[
f_{\la, \xi}(\eta)= f(\eta_\la\xi_{\la^c}).
\]
By definition \eqref{difi},  $\delta_i f_{\la, \xi}$ is non-decreasing when $\la$ grows, and is bounded by
$\delta_i f$.
By lemma \ref{boundlem}, it follows that $f_{\la,\xi}$ is bounded. 

 Therefore $\|\delta f_{\la,\xi}\|_2$ is bounded by $\|\delta f\|_2$ which is finite because $f\in \Delta^2(\Omega)$.
Therefore, using the assumed uniform convergence of $f_{\la,\xi}$ to $f$, and the assumed bound \eqref{GCB} for bounded local functions
in $\Delta^2(\Omega)$, we obtain (using dominated convergence)
\[
\int \e^{f-\int f \dd\mu}\dd\mu
= \lim_{\la\uparrow\Zd} \int \e^{f_{\la, \xi}-\int f_{\la,\xi} \dd\mu}\dd\mu
\leq \lim_{\la\uparrow\Zd}\e^{\tfrac{C}{2} \|\delta f_{\la, \xi}\|^2_2}
= \e^{\tfrac{C}{2} \|\delta f\|^2_2}.
\]
\end{proof}


\subsection{Relative entropy}

For a probability measure $\mu$, we denote by $\mu_\la$ its restriction to the sub-$\si$-algebra $\frakF_\la= \si\{ \eta_i, i\in \la\}$, generated by the projection $p_{\la}:\Omega\to\Omega_\la$.
We also denote by $\boB_\la$ the set of bounded $\frakF_\la$-measurable functions from $\Omega$ to $\R$.

For two probability measures $\mu,\mu'$ on $\Omega$ and $\la\in\boS$, we define the relative entropy
of $\mu'$ with respect to $\mu$ by
\[
\re_\la (\mu'|\mu)=
\begin{cases}
\mathlarger{\int} \dd\mu'_\la \log\myfrac[3pt][2pt]{\dd\mu'_\la}{\dd\mu_\la} & \text{if} \quad \mu'_\la\ll\mu_\la\\
+\infty & \text{otherwise}.
\end{cases}
\]
We further denote by $(\la_n)_{n\in\N}$ the sequence of cubes $\la_n=[-n,n]^d\cap\Zd$.
\bd[Lower relative entropy density]
For two probability measures $\mu,\mu'$ on $\Omega$, define the {\em lower relative entropy density} by
\[
\re_*(\mu'|\mu)=\liminf_n\frac{\re_{\la_n} (\mu'|\mu)}{|\la_n|} .
\]
\ed
We have the following variational characterization of relative entropy (for a proof see for instance \cite[p. 100]{blm})
\be
\label{vp-re}
\re_\la(\mu'|\mu)=\sup_{f}\left(\,\int f\dd\mu'_\la -\log \int \e^{f}\dd\mu_\la\right)
\ee
where the supremum is taken over all $\frakF_{\la}$-measurable functions such that $ \int \e^{f}\dd\mu_\la<\infty$.


\subsection{Main result}

In the main theorem below we prove that the Gaussian concentration bound implies strict positivity of the lower relative entropy density.

\bt\label{mainthm}
If $\mu$ is translation invariant and satisfies $\gcb{C}$ then for all $\mu'$ translation invariant, and $\mu'\not=\mu$ we have
\[
\re_*(\mu'|\mu) >0.
\]
\et
We start with a lemma from \cite{ccr}. For the reader's convenience, we repeat the short proof here.

\bl\label{young}
For $f$ such that $\|\delta f\|_1<+\infty$ and $\Lambda\in\boS$, we have
\[
\left\|\delta \left(\sum_{i\in \la}\tau_i f\right)\right\|_2^2\leq |\la|{\|\delta f\|_1^2}.
\]
\el
\begin{proof}
For $\la\subset\Zd$, let $\1_\la$ denote the indicator function of $\la$ (that is, $\1_\la(i)=1$ if $i\in\la$ and $\1_\la(i)=0$ otherwise).
Then for every $j\in\Zd$ we have
\[
0\leq \left[\delta \left(\sum_{i\in \la}\tau_i f\right)\right]_j \leq  \sum_{i\in\Zd} (\delta_{i+j} f) \1_\la (i)=(\delta f* \1_\la)_j.
\]
As a consequence, using Young's inequality for convolutions, we obtain
\[
\|\delta f\|_2^2 \leq \|\delta f* \1_\la\|^2_2\leq \| \1_\la\|_2^2\, \|\delta f\|_1^2= |\la| \|\delta f\|_1^2.
\]
\end{proof}

\begin{sproof}{of {\bf{\small{T}{\scriptsize HEOREM}}} \ref{mainthm}}
For the cube $\la_n$ and a bounded local function $f$ whose dependence set is included in the cube $\la_r$, for some $r$, it follows from \eqref{vp-re} that
\[
\frac{\re_{\la_{n+r} }(\mu'|\mu)}{|\la_n|}
\geq \frac{1}{|\la_n|} \left( \int \sum_{i\in\la_n} \tau_i f \dd\mu'
-\log \int \e^{\sum_{i\in\la_n} \tau_i f}\dd\mu\right)
\]
where we used that $\sum_{i\in\la_n} \tau_i f$ is measurable with respect to $\frakF_{\la_{n+r}}$.
Now if $\mu$ satisfies $\gcb{C}$ and both $\mu'$ and $\mu$ are translation invariant then we can estimate further as follows.
Start by noticing that, by combination of the assumed $\gcb{C}$ and Lemma \ref{young}, we have the estimate
\[
\log \int \e^{\sum_{i\in\la_n} \left(\tau_i f-\int f \dd\mu\right)}\dd\mu\leq \frac{C}{2}|\la_n| \|\delta f\|_1^2.
\]
As a consequence, using translation invariance of both $\mu$ and $\mu'$ we obtain
\begin{align}
\nonumber
& \frac{\re_{\la_{n+r} }(\mu'|\mu)}{|\la_n|} \\
\nonumber
& \geq  \frac{1}{|\la_n|} \left( \int \sum_{i\in\la_n} \tau_i f \dd\mu'
-\log \int \e^{\sum_{i\in\la_n} \tau_i f}\dd\mu\right)
\nonumber\\
&=\frac{1}{|\la_n|} \left( \int \sum_{i\in\la_n} \tau_i f \dd\mu' -\int \sum_{i\in\la_n} \tau_i f \dd\mu
-\log \int \e^{\sum_{i\in\la_n} (\tau_i f-\int \tau_i f \dd\mu)}\dd\mu\right)
\nonumber\\
& \geq \int f \dd\mu'-\int f\dd\mu -\frac{C}{2} \|\delta f\|_1^2.
\label{inter}
\end{align}
Consider a bounded local function $f$ such that $\int f \dd\mu'-\int f \dd\mu>u>0$ (this function
exists by the assumption that bounded local functions are measure separating). Put $\|\delta f\|_1^2=:\varrho$.
(Observe that $\varrho <\infty$ by assumption, and $\varrho\neq 0$ since $f$ cannot be a constant.) Assume that the dependence set of $f$ is included in the cube $\la_r$.
Replace $f$ by $\beta f$ in the inequality \eqref{inter}, and optimize over $\beta$. Then we obtain, for all $n\in \N$, the inequality
\begin{align}
\frac{\re_{\la_{n+r} }(\mu'|\mu)}{|\la_n|}
& \geq  \sup_{\beta\geq 0}\left\{\beta\left(\int f \dd\mu'-\int f\dd\mu\right) -C \beta^2\|\delta f\|_1^2\right\}
\nonumber\\
& \geq \sup_{\beta\geq 0} \left(\beta u- \frac{C}{2}\beta^2\varrho\right)= \frac{u^2}{2\,C\varrho}>0.
\end{align}
Since $r$ is fixed, we can take the limit inferior in $n$, and using  $|\la_n|/|\la_{n+r}|\to 1$ as $n\to\infty$, and we obtain

\[
\liminf_{n\to\infty}\frac{\re_{\la_{n+r} }(\mu'|\mu)}{|\la_n|} >0.
\]
\end{sproof}
From the proof of Theorem \ref{mainthm} we can infer the following more quantitative lower bound on the relative entropy density. Before stating, we define a distance between probability measures.
\bd
Define the following distance between probability measures
\be\label{dist}
d(\mu, \mu')=\sup \left\{ \int f\dd\mu -\int f\dd\mu' : f\in \caL(\Omega), \|\delta f\|_1\leq 1\right\}.
\ee
\ed
The metric defined above generates the quasilocal topology, and therefore convergence in this metric implies weak convergence.

Then we have the following quantitative version of the main theorem.
\bp
If $\mu$ satisfies $\gcb{C}$ then, for $\mu'$ translation invariant we have
\be\label{qunat}
\re_*(\mu'|\mu) \geq \frac{d(\mu', \mu)^2}{2C}.
\ee
\ep
\begin{proof}
From the proof of Theorem \ref{mainthm} we infer that for $f$ such that
\[
\int f \dd\mu'-\int f\dd\mu \geq u, \ \text{and}\ \|\delta f\|_1^2\leq \varrho
\]
we have
\[
\re_*(\mu'|\mu)\geq \frac{1}{2C} \left(\frac{u}{\sqrt{\varrho}}\right)^2.
\]
Therefore, for
\[
\int f \dd\mu'-\int f\dd\mu \geq \epsi \quad \text{and}\quad \|\delta f\|_1^2\leq 1
\]
we have
\be\label{sineq}
\re_*(\mu'|\mu)\geq \frac{\epsi^2}{2C}.
\ee
By definition of the distance \eqref{dist}, this is equivalent with the statement that $d(\mu',\mu)\geq \epsi$ implies
\eqref{sineq}.
This implies \eqref{qunat}.
\end{proof}

\br
Our distance $d(\mu,\mu')$ between probability measures ressembles the so-called Dobrushin distance, denoted by $D(\mu,\mu')$, which consists in taking the supremum of
$\int f \dd\mu-\int f\dd\mu'$ over a wider set of functions, namely $f$ is required to be measurable and such that $\|\delta f\|_1\leq 1$. Hence $d(\mu,\mu')\leq D(\mu,\mu')$ for a general
pair $\mu,\mu'$ of probability measures.  In the special case of finite $S$, one has $d=D$. In \cite{AG-M}, it is proved that $D$ is equal to what the authors called the Steiff distance $\overline{d}$ which is defined in terms of couplings, and which
generalizes the Ornstein distance. The equality between $D$ and $\overline{d}$ is reminiscent of the Kantorovich-Rubinstein duality theorem.
\er

\br
Inequality \eqref{qunat} is reminiscent of a well-known abstract inequality relating the relative entropy and the Wasserstein distance due to Bobkov and G\"otze \cite{bg}.
But our context is different, because we consider the thermodynamic limit and the relative entropy density.
Nevertheless, as shown in \cite{ccr}, we can exploit Bobkov-G\"otze theorem in the special case of finite $S$, putting the Hamming distance on $S^{\la_n}$, to get
\[
\re_*(\nu|\mu)\geq \frac{\bar{d}^{\,2}(\mu,\nu)}{2C}
\]
where 
\[
\bar{d}(\mu,\nu):=\lim_{n\to+\infty}  \frac{W_1(\mu_{\la_n},\nu_{\la_n};\db_n)}{|\la_n|}
\]
and $\db_n(\omega,\eta)=\sum_{i\in \la_n} \1_{\{\omega_i\neq\eta_i\}}$, $W_1(\mu_{\la_n},\nu_{\la_n};\db_n)$ being the Wassertein distance between $\mu_{\la_n}$ and $\nu_{\la_n}$.
\er

The following corollary is useful in the context of stochastic dynamics, where one can often show convergence of relative entropy density.
See e.g.\ \cite{kunch} and references therein.
\bc\label{convcol}
Let $\mu$ be a translation-invariant probability measure which satisfies $\gcb{C}$.
If $(\mu_n)$ is a sequence of translation invariant probability measures such that
\be\label{reldist}
\lim_{n\to\infty} \re_* (\mu_n|\mu)=0.
\ee
Then $\mu_n\to\mu$ in the sense of the distance \eqref{dist}, and therefore also $\mu_n\to\mu$ weakly. 
\ec
\bpr
By \eqref{qunat}, \eqref{reldist} implies
\[
\lim_{n\to\infty}d(\mu_n, \mu)^2\leq  2C \lim_{n\to\infty}\re_* (\mu_n|\mu)=0.
\]
Therefore we have convergence in the metric $d$, which, as we remarked before, implies weak convergence.
\epr

\section{Applications: uniqueness of equilibrium states and beyond}

In this subsection we provide some settings where we can conclude uniqueness of a set of ``(generalized) translation-invariant Gibbs measures''
via Theorem \ref{mainthm}. We start with the set of translation-invariant Gibbs measures associated to an absolutely summable potential.
Then we consider generalizations and modifications of such sets.

\subsection{Uniqueness of equilibrium states}
In this subsection we briefly introduce the necessary basics of Gibbs measures. The reader familiar with
the theory of Gibbs measure can skip this subsection.
The reader is referred to \cite{geo} (especially chapter 16) or \cite[Chapter 2]{enter} for more background on the Gibbs formalism.

Let $\lambda$ be a probability measure on $S$, and for $\lambda_\Lambda (\dd\si_\la)=\otimes_{i\in\Lambda} \lambda(\dd\si_i)$ the corresponding
product measure on $S^\la$.
The measure $\lambda$ is called the ``a priori'' measure on $S$, with associated a priori measure $\otimes_{i\in \Zd}\lambda (\dd\si_i)$ on $\Omega$.

We call a uniformly absolutely summable translation-invariant potential a function
\[
U: \boS\times \Omega\to \R
\]
such that
\bi
\item[(a)] Local potentials: for all $A\in \boS$, $U(A, \cdot)$ is  $\mathfrak{B}\!_A$-measurable and continuous.
\item[(b)] Absolutely summable: $\sum_{A\ni 0} \|U(A, \cdot)\|_\infty <\infty$.
\item[(c)] Translation invariant: $U(A+i, \tau_i \si)=U(A, \si)$ for all $\si\in \Omega$, $A\in \boS$.
\ei
Let us call $\caU$ the
set of uniformly absolutely summable translation invariant potentials.
Then we build the local Gibbs measures with boundary condition $\xi\in \Omega$. For a finite subset $\la\in \boS$, the Gibbs measure in volume
$\la$ with boundary condition $\xi$ outside $\la$ is defined via
\[
\gamma_\la(\dd\si_\la|\xi)= \frac{\e^{-H^\xi_\la(\si_\la)}}{Z_\la^\xi} \lambda_\la (\dd\si_\la)
\]
where $H^\xi_\la$ is the Hamiltonian in volume $\la$ with boundary condition $\xi$:
\[
H_\la^\xi(\si_\la)= \sum_{A\cap \la \not=\emptyset} U(A, \si_\la \xi_{\la^c})
\]
and where ${Z_\la^\xi}$ is the normalization
\[
{Z_\la^\xi}= \int{\e^{-H^\xi_\la(\si_\la)}}\lambda_\la (\dd\si_\la).
\]
The family $\gamma_\la(\dd\si_\la|\cdot), \la\in \boS$ is called the Gibbsian specification associated to
the potential $U$ (with a priori measure $\lambda$).

By the uniform absolute summability of $U$, we automatically have that for all
$f$ local and continuous,
the function $\xi\mapsto \int f(\si_\la\xi_{\la^c}) \gamma_\la (\dd\si_\la|\xi)$ is quasi-local and continuous.
We say that the specification $\gamma_\la(\cdot|\cdot)$ is quasi-local.

We then call a measure Gibbs $\mu$ with potential $U$ (and a priori measure $\lambda$  if $\gamma_\la(\dd\si_\la|\xi)$ if it
is consistent with the finite volume Gibbs measures, i.e., if for all $f:\Omega\to\R$ bounded and measurable,  and $\la\in \boS$
\[
\E_\mu (f|\mathfrak{B}_{\la^c})(\xi)= \int f(\si_\la \xi_{\la^c}) \gamma_\la (\dd\si_\la|\xi)
\]
$\mu$-almost every $\xi$.

We denote by $\caG_{\tau}(U)$ the set of translation invariant Gibbs measures associated to the potential $U$. These measures are called the ``equilibrium states'' associated to $U$.

The variational principle (\cite{geo} chapter 16) implies that if $\mu, \nu\in\caG_\tau(U)$ then
$s_*(\mu|\nu)=s_*(\nu|\mu)=0$, and conversely if $\mu\in \caG_\tau(U)$ and $\nu\in \caP_\tau$ is such that
$s_*(\nu|\mu)=0$ then $\nu\in \caG_\tau(U)$.
As a consequence of Theorem \ref{mainthm}, we then obtain the following:
\bp
Let $U\in \caU$.
If $\mu\in \caG_\tau(U)$ satisfies $\gcb{C}$ for some $C>0$ then
$\caG_\tau(U)=\{\mu\}$.
\ep
This substantially extends the implication between $\gcb{C}$ and uniqueness of equilibrium states
from \cite{moles} where we only considered finite single spin spaces $S$.

\subsection{Sets of zero information distance}

The example of the set of equilibrium states from the previous subsection leads naturally to the more general notion of ``zero information distance sets'' defined below.
\bd
We
call a subset $\caK\subset\caP_{\tau}(\Omega)$ a zero information distance set if
for all $\mu, \mu'\in \caK$, $s_*(\mu|\mu')=s_*(\mu'|\mu)=0$.
\ed
From Theorem \ref{mainthm} we then obtain immediately the following proposition.
\bp\label{zero}
Let $\caK\subset\caP_{\tau}(\Omega)$ be a zero information distance set. If there exists $\mu\in \caK$ which
satisfies $\gcb{C}$ for some $C\in (0,\infty)$, then $\caK$ is a singleton.
\ep
We now give three further examples (beyond equilibrium states) of such zero information distance sets, illustrating
Proposition \ref{zero}.

\begin{itemize}

\item[a)] Asymptotically decoupled measures and $\Pi^f$-compatible measures.
A first generalization of the Gibbsian context is provided in the realm of ``asymptotically decoupled measures'' via the notion of $\Pi^f$-compatible measures, see
\cite{pfister}. This setting  goes beyond quasi-local specifications and therefore includes many relevant examples of non-Gibbsian measures.

In this setting  the set of  $\Pi^f$-compatible measures  (associated to a local function $f$) is a zero information set
((see \cite{pfister} Theorem 4.1), and therefore, if this set contains an element $\mu$ satisfying
$\gcb{C}$, then it coincides with the singleton $\{ \mu\}$.

\item[b)] Renormalization group transformations of Gibbs measures.

Another important class of examples is the following. We say that a transformation $T:\caP_{\tau}(\Omega)\to\caP_{\tau}(\Omega')$ preserves
zero information distance sets if a zero information distance set is mapped by $T$ onto  a zero information distance set.
Important examples of such transformations $T$ are local and translation invariant renormalization group transformations studied in \cite{enter}, section 3.1, p 960, conditions T1-T2-T3. Examples of such transformations include block-spin averaging, decimation and stochastic transformations such as the Kadanoff transformation. Because the transformations are ``local and translation invariant probability kernels'', one immediately infers the property $s_*(\mu T|\nu T)\leq s_*(\mu| \nu)$.

In this setting,
Proposition \ref{zero} implies that if $U\in \caU$,  $\mu\in \caG_\tau(U)$ is an associated  translation-invariant Gibbs measure, and $\mu T$ satisfies $\gcb{C}$ for some $C\in (0,\infty)$, then $\nu=\mu T$ for all $\nu$ such that $s_*(\nu|\mu T)=0$. In particular, this implies that
$\mu'T=\mu T$ for all $\mu'\in \caG_\tau(U)$. Indeed, in that case
$s_*(\mu'T|\mu T)\leq  s_*(\mu'|\mu)=0$.

Notice that $\mu T$ can be non-Gibbs, therefore  the implication
$\nu=\mu T$ for all $\nu$ such that $s_*(\nu|\mu T)=0$ cannot be derived from the variational principle.

\item[c)] Stationary measures for Ising spin Glauber dynamics.

An additional example of a zero information distance sets is the set of stationary and translation-invariant measures for a (Ising spin, i.e.,
$S$ is finite) Glauber dynamics under the condition that this set contains at least one translation-invariant Gibbs measure as stationary measure, see \cite{kunch}, Section 4.
As a consequence of Proposition \ref{zero} we then conclude that if there exists a translation-invariant Gibbs measure $\nu$ as stationary measure, and
there exists a translation-invariant stationary measure $\mu$ satisfying $\gcb{C}$ for some $C>0$, then $\mu=\nu$ coincide, and $\mu$ is the unique
translation-invariant stationary measure.
Moreover, in this setting one can show that when starting the dynamics from a translation invariant initial measure $\mu$ and denoting
$\mu_t$ for the measure at time $t>0$, we have $\re_*(\mu_t|\nu)\to 0$ as $t\to\infty$.
Then, from corollary \ref{convcol} we obtain that $\mu_t\to \nu$ as $t\to\infty$ in the sense of the distance \eqref{dist}.
\end{itemize}

\section{Generalization}\label{gensection}

In the setting of Section \ref{buirbu}, without the additional assumption of finiteness of the diameter of $S$, the definition \eqref{difi} of the oscillation
of $f$ is no longer appropriate, because it becomes natural to include unbounded functions, which makes \eqref{difi} infinite, e.g. for $S=\R$, and $\Omega=S^{\Zd}$ equipped with a product of Gaussian measures, the function $f(\eta)=\eta_i$ should be a possible choice.
We consider now a general standard Borel $S$, which is such that for the product space
$\Omega=S^{\Zd}$, quasi-local bounded functions are measure separating.

In order to proceed, we therefore associate to a function $f:\Omega\to\R$ an abstract sequence of oscillations
$\delta f= \delta_i f, i\in\Zd$,  satisfying the following conditions.
\bd
We say that a map $\delta: \caC(\Omega)\to [0,\infty]^{\Zd}$ is an allowed sequence of oscillations if the following
four conditions are met.
\ben
\item Translation invariance: $(\delta(\tau_i f))_j=\delta_{i+j} f$.
\item Non-degeneracy: $\delta_i f$ is zero for a function $f$ if and only if $f$ does not depend on the $i$-th coordinate, i.e.,
$\delta_i(f)=0$ if and only if for all $\eta, \si$ such that $\eta_j=\si_j$ for all $j\not=i$, $f(\eta)=f(\si)$.
\item Monotonicity:
for $\xi\in\Omega$ and $f$ a bounded quasi-local function, we consider the local approximation of $f$ given by
\[
f_{\la, \xi} (\eta)= f(\eta_\la\xi_{\la^c})
\]
Then we require that for all $\xi$ and for all $\la$, $\delta_i f_\la\leq \delta_i f$.
\item Degree one homogeneity: $\delta_i (\beta f)=|\beta| \delta_i f$ for all $\beta\in\R$.
\een
\ed
Notice that condition 3 implies  that for given $f$, $\xi$, $\la\subset\la'$, and $i\in\Zd$ we have
 $\delta_i f_{\la, \xi}\leq \delta_i f_{\la',\xi}$. Indeed
notice that $(f_{\la',\xi})_{\la,\xi}= f_{\la, \xi}$ for $\la\subset\la'$.

The most natural examples different from \eqref{difi} is
\[
\delta_i f= \sup_{\substack{\si, \eta\,\in\,\Omega: \\ \si_j=\eta_j, \forall j\not= i}} \frac{f(\si)-f(\eta)}{d_S(\si_i, \eta_i)}.
\]
More generally one can define
\[
\delta_i f= \sup_{\substack{\si, \eta\,\in\,\Omega: \\ \si_j=\eta_j, \forall j\not= i}} \frac{f(\si)-f(\eta)}{\psi(\si_i, \eta_i)}
\]
where $\psi:S\times S\to [0,\infty)$ satisfies $\psi(s,s')=0$ iff $s=s'$.

For a given sequence of oscillations $\delta$, we call a function $\delta$-Lipschitz if
$\sup_i \delta_i f<\infty$.
We then introduce
\[
\Delta^p(\Omega)= \{ f\in \caC(\Omega): \|\delta f\|_p<\infty\}.
\]

\bd
Let $\Omega=S^{\Zd}$ where $S$ is a standard Borel space. Assume an allowed  sequence of oscillations $\delta $ is given.
Let $\mu$ be a probability measure on $\Omega$.
We say that $\mu$ the Gaussian concentration bound w.r.t. $\delta$ with constant $C>0$ (still abbreviated $\gcb{C}$), if for all bounded local functions we have
\[
\int \e^{f-\int f \dd\mu}\dd\mu\leq \e^{\tfrac{C}{2} \|\delta f\|^2_2}.
\]
\ed

We then have the following analogue of Theorem \ref{mainthm}. Because the proof follows exactly the same steps as the proof of Theorem \ref{mainthm}, we leave it to the reader.
\bt
Assume $\delta$ is an allowed vector of oscillations. Assume that the set of
bounded local $\delta$-Lipschitz functions is measure separating.
If $\mu$ is translation invariant and satisfies $\gcb{C}$ then for all $\mu'$ translation invariant, and $\mu'\not=\mu$ we have
\[
\re_*(\mu'|\mu) >0.
\]
\et

As a final comment, we remark that the fact we have chosen the group $\Zd$ is for the sake of simplicity. We can work with more general amenable groups as in \cite{tempelman}.


\end{document}